\documentclass[12pt]{article}
\usepackage[cp1251]{inputenc}
\usepackage[russian]{babel}
\usepackage{amssymb}
\usepackage{amsmath}

\begin{document}

\begin{center}
\textbf{\large On $l^p$ -multipliers of functions analytic in the
disk} \footnote{This study was carried out within The National
Research University Higher School of Economics' Academic Fund
Program in 2013-2014, research grant No. 12-01-0079.}
\end{center}

\begin{center}
Vladimir Lebedev
\end{center}

\begin{quotation}
{\small \textbf{Abstract.} We consider bounded analytic
functions in domains generated by sets that have
Littlewood--Paley property. We show that each such function is
an $l^p$ -multiplier.

  References: 12 items.

  \textbf{Key words:} bounded analytic functions,
$l^p$ -multipliers, Littlewood--Paley property.

MCS 2010: Primary 42A45; Secondary 30H05.}
\end{quotation}

\quad

  Given a function $f$ analytic in the unit disk
$D=\{z\in\mathbb C : |z|<1\}$ of the complex plane $\mathbb C$,
consider its Taylor expansion:
$$
f(z)=\sum_{n\geq 0} \widehat{f}(n) z^n, \qquad z\in D.
\eqno(1)
$$
For $1\leq p\leq\infty$ let $A_p^+(D)$ denote the space of all
functions (1) such that the sequence of Taylor coefficients
$\widehat{f}=\{\widehat{f}(n), ~n=0, 1, \ldots\}$ belongs to
$l^p$. For $f\in A_p^+(D)$ we put
$\|f\|_{A_p^+(D)}=\|\widehat{f}\|_{l^p}$. A function $m$ analytic
in $D$ is called an $l^p$ -multiplier if for every function $f$ in
$A_p^+(D)$ we have $m\cdot f\in A_p^+(D)$. We denote the class of
all these multipliers by $M_p^+(D)$. This class is a Banach
algebra with respect to the natural norm
$$
\|m\|_{M_p^+(D)}=\sup_{\|f\|_{A_p^+(D)}\le 1} \|m\cdot
f\|_{A_p^+(D)}
$$
and the usual multiplication of functions. The classes $M_p^+(D)$
were studied in [1]--[6]. \footnote{There is a minor inconsistency
in \S~6 of the author's work [6]. Instead of the written ``the
Poisson integral'' there should be ``the Riesz projection
$\sum_{k=-\infty}^\infty c_k e^{ikt}\rightarrow\sum_{k\geq 0}
c_kz^k$~''.} We note that the case when $p\neq 1, \infty, 2$ is of
a special interest. It is well known that $M_p^+(D)=M_q^+(D)$ if
$1/p+1/q=1$, and
$$
A_1^+(D)=M_1^+(D)=M_\infty^+(D)\subseteq M_p^+(D)
\subseteq M_2^+(D)=H^\infty(D),
$$
where $H^\infty(D)$ is the Hardy space of bounded analytic
functions in $D$.

  Let $\Omega\subseteq \mathbb C$ be an open domain which contains the
disk $D$. We shall present a class of domains $\Omega$ such that
each bounded analytic function in $\Omega$ belongs to $M_p^+(D)$.
The case when $\Omega$ contains the closure of $D$ is trivial; in
this case each bounded analytic function in $\Omega$ belongs to
$A_1^+(D)$ and hence belongs to $M_p^+(D)$ for all $p, ~1\leq
p\leq\infty$. The nontrivial case is the case when the boundary of
$\Omega$ has common points with the boundary $\partial
D=\{z\in\mathbb C : |z|=1\}$ of the disk $D$.

   It was shown by Vinogradov [2] that if $r>1, ~0\leq\alpha<\pi/2,$
and $m$ is a bounded analytic function in the domain
$$
\Omega_0=\{z\in\mathbb C: |z|<r, ~\alpha<\arg (z-1)<2\pi -\alpha\},
\eqno(2)
$$
then $m\in\bigcap_{1<p<\infty} M_p^+(D)$. Using this result
Vinogradov gave the first examples of nontrivial (i.e., infinite)
Blaschke products in $M_p^+(D)$. Note that the boundary of a
domain (2) has only one common point with the boudary of $D$
(namely the point $z=1$). As we shall see, a statement similar to
Vinogradov's result holds for domains of a much wider class.
Functions analytic in the domains considered below can have
uncountable set of singularities on the boundary of the disk. \footnote{We note by the way
that the condition $\alpha<\pi/2$ in the Vinogradov theorem is
essential. For example, the function $S(z)=\exp \{(z+1)/(z-1)\}$
is bounded in the half-plane $\{z\in\mathbb C: \textrm{Re}
~z<1\}$, but as it was shown by Verbitski\v{\i} [4] $S$ belongs to
$M_p^+(D)$ only in the trivial case $p=2$.}

  As is usual, for an arbitrary domain $\Omega\subseteq\mathbb C$
by $H^\infty(\Omega)$ we denote the Hardy space of all bounded
analytic functions in $\Omega$. For $g\in H^\infty(\Omega)$ we put
$\|g\|_{H^\infty(\Omega)}=\sup_{z\in\Omega} |g(z)|$.

  Let $J$ be an arc in the boundary circle $\partial D$. Assume that
the length $|J|$ of $J$ is strictly less than $\pi$. Let $T_J$
be an arbitrary open isosceles triangle, whose base is the
chord that spans the arc $J$, and whose two other sides lie
outside $D$. By $\theta_{T_J}$ we denote the angle between
$\partial D$ and a side of $T_J$ distinct from its base.

   Consider an arbitrary closed set $F\subseteq \partial D$. Let
$\tau(F)$ be the family of all arcs complimentary to $F$ (i.e., of
all connected components of the compliment $\partial D\setminus
F$). We assume that each arc of the family $\tau(F)$ has length
strictly less than $\pi$. Consider the domain
$$
\Omega_F=D\cup\bigcup_{J\in\tau(F)} T_J,
$$
and require in addition that
$$
\inf_{J\in\tau(F)}\theta_{T_J}>0.
\eqno(3)
$$
We call a domain $\Omega_F$, obtained in this way, a star-like
domain generated by $F$.

  We shall show that under a certain condition imposed on
a set $F\subseteq\partial D$, every function, bounded and analytic
in $\Omega_F$, belongs to $M_p^+(D)$.

  Let $E$ be a closed set of Lebesque measure zero in the line
$\mathbb R$. Consider the family $\tau(E)$ of all intervals
complimentary to $E$ (i.e., of all connected components of the
compliment $\mathbb R\setminus E$). For an arbitrary interval
$I\subseteq\mathbb R$ define the operator $S_I$ by
$$
\widehat{S_I(f)}=1_I\widehat{f}, \qquad f\in L^p\cap L^2(\mathbb R),
$$
where $~\widehat{}~$ is the Fourier transform and $1_I$ is the
characteristic function of $I$ (i.e., $1_I(t)=1$ for $t\in I$,
$1_I(t)=0$ for $t\notin I$). Following [7], we say that a set $E$
has property $\mathrm{LP}(p) ~(1<p<\infty)$ if the corresponding
Littlewood--Paley quadratic function
$$
S(f)=\bigg (\sum_{I\in\tau(E)}|S_I(f)|^2\bigg )^{1/2}
$$
satisfies
$$
c_1(p)\|f\|_{L^p(\mathbb R)}\leq
\|S(f)\|_{L^p(\mathbb R)}\leq c_2(p)\|f\|_{L^p(\mathbb R)},
\qquad f\in L^p(\mathbb R)
$$
(with positive constants $c_1(p), ~c_2(p)$ independent of $f$). In
the case when $E$ has property $\mathrm{LP}(p)$ for all $p,
~1<p<\infty,$ we say that $E$ has property $\mathrm{LP}$.

   Let now $F$ be a closed set of measure zero in the boundary
circle $\partial D$. We say that $F$ has property $\mathrm{LP}(p)$
or property $\mathrm{LP}$ if $F=\{e^{it}, ~t\in E\}$, where
$E\subseteq [0, 2\pi]$ is a set that has property $\mathrm{LP}(p)$
or property $\mathrm{LP}$, respectively.

\quad

 \textsc{Remark 1.} A classical example of an infinite set
$E\subseteq\mathbb R$ that has property $\mathrm{LP}$ is $E=\{\pm
2^{k}, ~k\in\mathbb Z\}\cup\{0\}$, where $\mathbb Z$ is the set of
integers. At the same time there exist uncountable sets that have
property $\mathrm{LP}$. This was first established by Hare and
Klemes [8]. The existence of such sets was also noted in [9], see
details in [10, \S~4]. Let us state the corresponding result for
sets in $\partial D$. For each $p, ~1<p<\infty,$ there is a
constant $\beta_p ~(0<\beta_p<1)$ such that the following holds.
Let $F\subseteq\partial D$ be a closed set of measure zero.
Suppose that the arcs $J_k, ~k=1, 2, \ldots,$ complimentary to
$F$, being enumerated so that their lengthes do not increase,
satisfy $|J_{k+1}|/|J_k|\leq \beta_p$ for all sufficiently large
$k$. Then $F$ has property $\mathrm{LP}(p)$. This in turn implies
that if $\lim_{k\rightarrow\infty}|J_{k+1}|/|J_k|=0$, then $F$ has
property $\mathrm{LP}$.

\quad

  The result of this note is the following theorem.

\quad

\textsc{Theorem.} \emph{Suppose that a set $F\subseteq
\partial D$ has property $\mathrm{LP}(p)$, and $\Omega_F$ is
a star-like domain generated by $F$. Then $H^\infty
(\Omega_F)\subseteq M_p^+(D)$. If $F$ has property $\mathrm{LP}$,
then $H^\infty (\Omega_F)\subseteq \bigcap_{1<p<\infty}
M_p^+(D)$.}

\quad

\textsc{Proof.}  Let $G$ be an Abelian group and let $\Gamma$ be
the group dual to $G$. Consider a function $m\in L^\infty
(\Gamma)$ and the operator $Q$ defined by
$$
\widehat{Qf}=m\widehat{f}, \qquad f\in L^p\cap L^2 (G),
$$
where $~\widehat{}~$ stands for the Fourier transform on $G$. The
function $m$ is called an $L^p$ -Fourier multiplier if the
corresponding operator $Q$ is a bounded operator from $L^p(G)$ to
itself ($1\leq p\leq\infty$). Denote the class of all these
multipliers by $M_p(\Gamma)$ and put
$\|m\|_{M_p(\Gamma)}=\|Q\|_{L^p(G)\rightarrow L^p(G)}$. The
relation between the multipliers on the line $\mathbb R$ and on
the circle $\mathbb T=\mathbb R/2\pi\mathbb Z$ is well known [11]
(see also [12]). We shall need the Jodeit theorem [12] on the
periodic extension of multipliers. According to this theorem, if
$f\in M_p(\mathbb R)$ is a function that vanishes outside of the
interval $[0, 2\pi]$ and $g$ is the $2\pi$ -periodic function that
coincides with $f$ on $[0, 2\pi]$, then $g\in M_p(\mathbb T)$.
Note that there is a direct relation between the spaces $M_p^+(D)$
and $M_p(\mathbb T)$. Given a function $m\in H^\infty(D)$ consider
its (non-tangential) boundary function $m^*(t)=m(e^{it})$. The
conditions $m\in M_p^+(D)$ and $m^*\in M_p(\mathbb T)$ are
equivalent [3] (see also [5]).

  We shall also need the following statement. Let
$E\subseteq \mathbb R$ be a set that has property
$\mathrm{LP}(p)$. Suppose that a function $f\in L^\infty(\mathbb
R)$ is continuously differentiable on each interval complimentary
to $E$, and its derivative $f'$ satisfies
$$
|f'(t)|\leq \frac{c}{\textrm{dist}(t, E)},
\qquad t\in\mathbb R\setminus E,
\eqno(4)
$$
where $\textrm{dist}(t, E)$ stands for the distance from a point
$t$ to the set $E$ and $c>0$ does not depend on $t$. Then $f\in
M_p(\mathbb R)$.  This result of Sj\"ogren and Sj\"olin [7]
generalizes the well known Mikhlin--H\"ormander theorem.

   We note now that condition (3) implies the existance of
a constant $c=c(\Omega_F)>0$ such that if
$e^{it}\in\partial D\setminus F$, then the circle centered at
$e^{it}$ and of radius $r(t)=c\cdot \mathrm{dist}(e^{it}, F)$ lies
in $\Omega_F$. Denote this circle by $\gamma(t)$. Let $m\in
H^\infty(\Omega_F)$. Consider an arc $J$ complimentary to $F$. Let
$e^{it}\in J$. Consider the corresponding circle $\gamma(t)$. For
an arbitrary point $z$ that lies inside $\gamma(t)$ we have
$$
m'(z)=\frac{1}{2\pi i}\int_{\gamma(t)} \frac{m(\zeta)}{(\zeta -
z)^2}d\zeta.
$$
In particular,
$$
m'(e^{it})=\frac{1}{2\pi i}\int_{\gamma(t)} \frac{m(\zeta)}{(\zeta -
e^{it})^2}d\zeta.
$$
Hence, for the derivative $(m^*)'$ of the boundary function
$m^*(t)=m(e^{it})$ we obtain
$$
|(m^*)'(t)|=|ie^{it}m'(e^{it})|=\bigg|\frac{1}{2\pi i}\int_{\gamma(t)}
\frac{m(\zeta)}{(\zeta - e^{it})^2}d\zeta \bigg |\leq
$$
$$
\leq \frac{1}{2\pi}\int_{\gamma(t)} \frac{|m(\zeta)|}{|\zeta -
e^{it}|^2}|d\zeta|\leq \frac{1}{2\pi} 2\pi r(t)
\|m\|_{H^\infty(\Omega_F)}
\frac{1}{(r(t))^2}=
$$
$$
=c_1(\Omega_F)\|m\|_{H^\infty(\Omega_F)}
\frac{1}{\textrm{dist}(e^{it}, F )}.
\eqno(5)
$$

  Let $E\subseteq [0, 2\pi]$ be a set such that $F=\{e^{it}, ~t\in
E\}$ and $E$ has property $\mathrm{LP}(p)$. Without loss of
generality we can assume that $E$ contains the points $0$ and
$2\pi$. Define a function $f$ on $\mathbb R$ by $f(t)=1_{[0,
2\pi]}(t)m^*(t), ~t\in\mathbb R$. We see that (see (5)) the
function $f$ satisfies (4). Therefore, by the Sj\"ogren--Sj\"olin
theorem, we have $f\in M_p(\mathbb R)$. Hence, using the Jodeit
theorem, we obtain $m^*\in M_p(\mathbb T)$. Taking into account
the relation between multipliers on $\mathbb T$ and multipliers of
functions analytic in the disk $D$, we obtain $m\in M_p^+(D)$.

\quad

\textsc{Remark 2.} As far as the author knows, the question on the
existence of a set that has property $\mathrm{LP}(p)$ for some $p,
~p\neq 2,$ but does not have property $\mathrm{LP}$ is open.

\quad

\begin{center}
\textsc{References}
\end{center}

\flushleft
\begin{enumerate}

\item S. A. Vinogradov, ``Multiplicative properties of power
    series with coefficient sequence from $l^p$'', \emph{Dokl.
    Akad. Nauk SSSR}, \textbf{254}:6 (1980), 1301--1306.
    English transl. in \emph{Soviet Math. Dokl.},
    \textbf{22}:2 (1980).

\item S. A. Vinogradov, ``Multipliers of power series with the
    sequence of coefficients from $l^p$'', \emph{Zap. Nauchn.
    Sem. LOMI}, \textbf{39} (1974), 30-39; English transl. in
    \emph{J. Soviet Math.}, \textbf{8}:1 (1977).

\item N. K. Nikol'ski\v{\i}, ``On the spaces and algebras of
    Toeplitz matrices acting in $l^p$\,'', \emph{Sibirsk. Mat.
    Z.}, \textbf{7} (1966), 146--158; English transl. in
    \emph{Sibirian Math. J.}, 7 (1966).

\item I. E. Verbitski\v{\i}, ``Multipliers of spaces
    $l^p_A$'', \emph{Funkts. Anal. Prilozhen.}, \textbf{14}:3
    (1980), 67--68. English transl. in \emph{Functional
    Analysis and its Applications}, \textbf{14}:3 (1980),
    219--220.

\item V. V. Lebedev, ``Inner functions and $l^p$
    -multipliers'', \emph{Funkts. Anal. Prilozhen.},
    \textbf{32}:4 (1998), 10--21. English transl. in
    \emph{Functional analysis and its applications},
    \textbf{32}:4 (1998), 227--236.

\item V. V. Lebedev, ``Spectra of inner functions and
    $l^p$-multipliers'', \emph{Operator Theory: Advances and
    Applications}, \textbf{113} (2000), 205--212.

\item P. Sj\"ogren and P. Sj\"olin, ``Littlewood--Paley
    decompositions and Fourier multipliers with singularities
    on certain sets'', \emph{Ann. Inst. Fourier, Grenoble},
    \textbf{31}:1 (1981), 157--175.

\item K. Hare and I. Klemes, ``On permutations of lacunary
    intervals'', \emph{Trans. Amer. Math. Soc.},
    \textbf{347}:10 (1995), 4105--4127.

\item V. Lebedev and A. Olevski\v{\i}, ``Bounded groups of
    translation invariant operators'', \emph{C. R. Acad. Sci.
    Paris, Ser. I}, \textbf{322} (1996), 143--147.

\item V. V. Lebedev,  A. M. Olevski\v{\i}, ``$L^p$ -Fourier
    multipliers with bounded powers'', \emph{Izvestiya RAN:
    Ser. Mat.}, \textbf{70}:3 (2006), 129-166; English transl.
    in \emph{Izvestya: Mathematics}, \textbf{70}:3 (2006),
    549-585.

\item K. de Leew, ``On $L^p$-multipliers'', \emph{Ann. Math.},
    \textbf{81} (1965), 364--379.

\item M. Jodeit, ``Restrictions and extensions of Fourier
    multipliers'', \emph{Studia Mathematica}, \textbf{34}
    (1970), 215-226.

\end{enumerate}

\quad

\quad

\qquad National Research University Higher School of Economics,\\
\qquad Moscow Institute of Electronics and Mathematics

\qquad e-mail: \emph {lebedevhome@gmail.com}

\end{document}